\documentclass[12pt]{amsart}
\usepackage{fullpage}
\usepackage{amssymb}
\usepackage{mathrsfs}
\usepackage{amsmath}
\usepackage{hyperref}
\usepackage{amsthm}
\usepackage{amssymb}
\usepackage{color}
\usepackage[english]{babel}
\usepackage{indentfirst}

\newtheorem{theorem}{Theorem}[section]

\newtheorem{lemma}[theorem]{Lemma}
\newtheorem{proposition}[theorem]{Proposition}

\theoremstyle{definition}

\theoremstyle{remark}

\newcommand{\be}{\begin{equation}}
\newcommand{\ee}{\end{equation}}
\numberwithin{equation}{section}

\newcommand{\R}{\mathbb R}

\newcommand{\T}{\mathbb T}

\newcommand{\co}{\colon}
\newcommand{\ep}{\varepsilon}

\newcommand{\bbp}{\mathbf{\bar p}}
\newcommand{\bbq}{\mathbf{\bar q}}
\newcommand{\bbc}{\mathbf{\bar c}}
\newcommand{\bbz}{\mathbf{\bar 0}}

\title[Entropy non-expansive KAM systems]{Key ideas behind perturbing any completely integrable Hamiltonian system obtaining volume entropy non-expansiveness}
\author{Dmitri Burago, Dong Chen and Sergei Ivanov}

\address{Dmitri Burago: Pennsylvania State University, Department of Mathematics, University Park, PA 16802, USA}
\email{burago@math.psu.edu}
\address{Dong Chen: Pennsylvania State University, Department of Mathematics, University Park, PA 16802, USA}
\email{dxc360@psu.edu}
\address{Sergei Ivanov: St.Petersburg Department of Steklov Mathematical Institute, Russian Academy of Sciences, Fontanka 27, St.Petersburg 191023, Russia}
\email{svivanov@pdmi.ras.ru}
\date{}    

\begin{document}
\maketitle

\begin{abstract}

For an arbitrary completely integrable Hamiltonian system with $n$ degrees of freedom $(n \geq 2),$ we construct a $C^{\infty}$-small Hamiltonian perturbation with positive volume entropy.
The key point is that the perturbation can be made near any Liouville torus and that the positive entropy can be generated in any arbitrarily small tubular neighborhood of one trajectory.   We explain the key ideas of our construction, while the details of the proofs are to appear elsewhere.  The novelty of the paper is \textcolor{black}{developing tools for}
perturbing the flow in such a way that a given (symplectic) diffeomorphism becomes the Poincar\'{e} return map near a closed trajectory.  \textcolor{black}{The full version can be found in  \cite{BCI}, and we hope this paper may assist in reading it or just get the basic ideas. }
\end{abstract}

\section{Introduction}
\footnote{2010 \emph{Mathematics Subject Classification.} 37A35, 37J40, 53C60.

\emph{Key words and phrases.} volume entropy non-expansive maps, KAM theory, Hamiltonian flow, perturbation.

The first author was partially supported
by NSF grant DMS-1205597. The second author was partially supported by Dmitri Burago's Department research fund 42844-1001. The third author was partially supported by RFBR grant 20-01-00070.}



\textcolor{black}{A rather technical paper  \cite{BCI} addresses a long-standing problem in Hamiltonian dynamics. We want to stay away from technicalities and explain what is behind our argument. It is not the story of constructing a single symplectomorphism (which is $C^{\infty}$ close to identity) with positive volume entropy.  It is rather about embedding it into the Poincar\'{e} map of some Hamiltonian flow, and the construction of this embedding relies on a generalization of techniques in  \cite{BI}. }

\textcolor{black}{We start with the most ``deterministic" Hamiltonian systems (sometimes called ``stable"  in other literature, for example,  \cite{M2}).
A Hamiltonian flow on a symplectic manifold $(\Omega^{2n},\omega_0)$} is said to be \textit{completely integrable} if it enjoys $n$ 
algebraically independent first integrals which pairwisely Poisson commute.   Here  ``algebraically independent" means that the gradient vectors are linear independent almost everywhere,  \textcolor{black}{ while ``Poisson commute" stands for the commutativity of the induced Hamiltonian flows}.
\textcolor{black}{These $n$ pairwisely commuting Hamiltonian flows define an $\mathbb{R}^n$-action on the phase space. This action is proper due to algebraic independence, and its compact orbits are in fact tori (called the  \textit{Liouville Tori}),  according to the Liouville-Arnold Theorem (a precise statement can be found in \cite[Theorem 5.5.21]{KH}). }
The motion on each Liouville tori is conjugate to a linear flow on a standard torus. These invariant tori are in fact the common level sets of the angle variables in the action-angle coordinates.

\textcolor{black}{What happens if we make a tiny perturbation on these ``deterministic"  systems? At first glance, it seems unlikely those invariant tori would persist under even $C^{\infty}$-small perturbation.  For example, it is not hard to break the periodic Liouville tori by slightly changing the Hamiltonian. However,  the celebrated Kolmogorov-Arnold-Moser(KAM) Theorem \cite{Ar1}\cite{K}\cite{M} shows that any $C^{\infty}$-small perturbation of any nondegenerate (in a sense which is not relevant to this brief exposition) completely integrable system results in an overwhelming measure of invariant tori (called \textit{KAM tori}) on which the dynamics 
 is still quasi-periodic.}
``An overwhelming measure" means that the measure of the tori which do not survive goes to zero as the size of the perturbation diminishes. The tori which survive have rotation vectors whose directions are ``sufficiently irrational" (a certain degree of being {\it Diophantine}, the precise condition is a bit technical and of no relevance here).  
 
The KAM tori can be viewed as the ``deterministic" part   of the phase space in the classical model where the error in initial conditions does not grow exponentially.  
\textcolor{black}{This paper is concerned with the ``chaotic" ( or ``unstable") part.  In fact, it is known that some chaotic behavior can appear outside the KAM tori.  More precisely,  a $C^2$-generic Hamiltonian flow has positive topological entropy \cite{N1} (such genericity were later confirmed in  \cite{KW}  and  \cite{C} for geodesic flows).   To get positive topological entropy, however,  it suffices to find one Smale horseshoe (even of zero measure with respect to the canonical invariant measure),  
for which we only need one transverse homoclinic point.
}

\textcolor{black}{ 
A hyperbolic set with positive measure leads to  coexistence of two totally different behaviors:  integrable (or ``laminary") and ``strongly chaotic" (or ``turbulent").  The latter property can be characterized by positive volume entropy.
The \textit{volume entropy}, also known as \textit{metric entropy} or \textit{Kolmogorov-Sinai entropy},  is the measure-theoretic 
entropy with respect to the Liouville measure on a level set (or with respect to the symplectic volume in the entire phase space).} \textcolor{black}{ The topological entropy measures the ``maximal exponential complexity", while the volume entropy measures the  ``average exponential complexity" (a precise definition of both entropies can be found in Chapters 3 and 4 in \cite{KH}).  It is a way more difficult to generate positive volume entropy (positivity of the volume entropy implies that for the topological entropy but not vice versa  \cite{BT1}).  It had been open for decades whether the perturbed systems in KAM Theorem may admit positive volume entropy,  until recently in \cite{Ch} the second author managed to generate positive volume entropy by perturbing the standard metrics on $\mathbb{T}^n (n\geq 3)$, using the techniques developed in \cite{BI}.}

\textcolor{black}{ Here we apply a generalization of the techniques in \cite{BI}  to Lagrangian submanifolds to   get positive volume entropy by perturbing the flow near any Liouville torus:}

\begin{theorem}\label{thm2}
Let $\Phi^t_H$ be a completely integrable Hamiltonian flow on 
a symplectic manifold $\Omega=(\Omega^{2n},\omega)$ with $n\geq 2$. 
For any Liouville torus $\mathcal T\subset\Omega$,
one can find a $C^{\infty}$-small perturbation $\widetilde{H}$ of $H$ 
such that the resulting Hamiltonian flow $\Phi^t_{\widetilde{H}}$ 
has positive volume entropy.
Furthermore, such perturbation can be made in an arbitrarily small neighborhood
of $\mathcal T$ and such that the flow is entropy non-expansive. If $\Omega$ is compact,
there are perturbations which produce positive volume entropy with respect to $\omega^n$ on the whole $\Omega$.
\end{theorem}

{\it Remark}. Theorem \ref{thm2} answers perhaps the most intriguing question in the KAM Theory. The question, which is probably due to A. Kolmogorov, asks if arbitrarily small perturbations of a completely integrable system
may result in dynamics with positive volume entropy.  As one of the referees points out,    what Kolmogorov considered at the 50s might be real-analytic perturbations. In this case, our method does not work \textcolor{black}{since we use bump functions.  We do not know whether we can get positive volume entropy from real-analytic perturbations and, even more intriguing,  for  generic perturbations.}

{\it Remark}.  \textcolor{black}{What is very important to us, our examples are entropy non-expansive.} This resolves another well-known problem.  One says a flow $\Phi^t$ is \textit{entropy non-expansive} if the positive volume entropy can be generated in an arbitrarily small tubular neighborhood of one orbit.  
\textcolor{black}{This notion was introduced by Bowen \cite{B} (and by Knieper \cite{K98} for geodesic flows),  and independently by the first author \cite{Bu}. }
Entropy non-expansiveness is a bit counter-intuitive since hyperbolic dynamics tends to expand and occupy all space.   \textcolor{black}{In fact, when Bowen defined entropy expansiveness in \cite{B}, he confessed he knew no example of entropy non-expansive diffeomorphism of a compact manifold.} \textcolor{black}{In our example,  however,   the positive volume entropy is generated even  near a periodic orbit. }

{\it Remark}. 
 \textcolor{black}{When the unperturbed Hamiltonian (Lagrangian) flows are geodesic flows on the cotangent bundle of flat Finsler tori, one can also make perturbations in the class of Finsler metrics to get positive volume entropy. However, we do not know whether the resulting metrics can be made Riemannian if the unperturbed geodesic flows are Riemannian.  This remains a mystery. }

There are significant differences between examples in \cite{BI}  \cite{Ch} and the one in this paper. 
Papers  \cite{BI} and  \cite{Ch} only deal with very specific systems: the geodesic flows of the standard metrics on 
$\mathbb S^n$ and $\mathbb T^n$;
the flow on standard $\mathbb S^n$
is KAM-degenerate (and periodic), and the example on $\T^n$ is entropy expansive.  In contrast to that, here
we can work with any completely integrable system, regardless of whether
it is KAM-nondegenerate or not,
and the perturbed flow is entropy non-expansive.

Also, thanks to the  results in Berger-Turaev \cite{BT2}, our result applies to $n=2$.  In this case,  the $2$-dimensional KAM tori separate the $3$-dimensional energy level. Thus 
no Arnold diffusion is possible in such systems.  The existence of positive volume entropy between these tori is a little bit surprising.   \textcolor{black}{Moreover,   the result in \cite{BT2} can be generalized to perturbing identity maps on higher dimensional discs, allowing us to address the cases when $n\geq 2$}.  Before \cite{BT2},   \textcolor{black}{for $n\geq 4$, }we perturbed the identity map by symplectically embedding the geodesic flow in the cotangent bundle of a negatively curved surface into Euclidean space (see  \cite[Lemma 5.1]{BI}), which inevitably required larger dimension.  One of the features of this approach,
however, is that we have more flexibility of perturbations than those in  \cite{BT2}. Of course, 
there is an obvious infinite-dimensional space of perturbations obtained by conjugating any example by
a symplectomorphism; here we can, however, make {\em essential} changes by varying the metric of the negatively curved surface we employ.

The above improvements run into serious difficulties and required new techniques and ideas.  The first challenge is how to generate positive volume entropy by perturbing the Poincar\'{e} map $R$.   \textcolor{black}{The result in \cite{BT2}  is not enough to guarantee positive entropy on various level sets since the volume entropy may vanish under perturbations.}
A more serious challenge is to realize the perturbed return map $\tilde{R}$ as the Poincar\'{e} map of some Hamiltonian $\tilde{H}$ that is $C^{\infty}$-close to  $H$. In the special case when $\Phi^t_H$ is the geodesic flow on a Finsler manifold, the \textcolor{black}{techniques} in \cite{BI} can be applied to get the desired perturbation $\tilde{H}$.  
\textcolor{black}{If the return map $\tilde{R}$ is the time-one map of some flow,  one could construct $\tilde{H}$ using the suspension of the generating vector field (see e.g. \cite{CHPZ21}).  However,  the map in Berger-Turaev \cite{BT2}, as well as the $\tilde{R}$ in our construction, is not generated by a flow. Thus we introduce new tools to build up $\tilde{H}$.} 

The argument consists of  two parts outlined below:
\vspace*{-1.5mm}
\begin{itemize}
\item In Section \ref{sec: proof} we perturb the Poincar\'{e} map to create a ``\textit{periodic spot}" (a disc consisting of points with the same period, see say \cite{GT10}\cite{GT17}\cite{T10}) in the cross-section.  We overcome the first challenge above by building up  special symplectic coordinates on each level sets and then apply the Morse-Bott Lemma (see \cite{BH} for the precise statement) to guarantee the perturbations on different level sets behave in the same way.  We then insert a positive-entropy symplectomorphism of the disc in the periodic spot using the results in  \cite{BT2}.  After all, the original Poincar\'{e} map is perturbed to get positive volume entropy and the resulting map is entropy non-expansive. 

\item We extend the perturbed Poincar\'{e} map to a perturbation of the Hamiltonian in Section \ref{sec:R to H} and address the second challenge.  The major novelty in this part is using Lagrangian submanifolds \textcolor{black}{to perturb} a Hamiltonian flow so that \textcolor{black}{a given (symplectic) diffeomorphism becomes the Poincar\'{e}  return maps of some Hamiltonian flow near a periodic trajectory}
 (see Proposition \ref{perturb all levels}). 
\end{itemize}

\textcolor{black}{The volume entropy we get in Theorem \ref{thm2} is inevitably small since  the measure of the ``chaotic" part is small and volume entropy is the ``average of exponential instability".  Nevertheless it would be interesting to explore how large entropy can be generated depending on the size of perturbation (any estimate would certainly involve some characteristics of the unperturbed system)? Probably some (very non-sharp) lower bounds can be obtained by a careful analysis of the proof.  }


\textit{Acknowledgments}.
We are grateful to Leonid Polterovich for his valuable help for boosting our understanding of the technique related to Lagrangian submanifolds. We also thank Vadim Kaloshin, late Anatole Katok, Federico Rodriguez Hertz,  Yakov Sinai,  Pierre Berger,  and Dmitry Turaev for useful discussions.  We would also like to thank the anonymous referees for many suggestions for improving the manuscript.

\section{Hamiltonian perturbations with prescribed Poincar\'e maps}
\label{sec:R to H}

In this section we show that certain perturbations of Poincar\'e maps can be realized as  Poincar\'e maps  of perturbed Hamiltonian flows.  Let
$\Omega=(\Omega^{2n},\omega)$ be a symplectic manifold, $n\ge 2$,
$H\co\Omega\to\R$ a Hamiltonian, and $\{\Phi_H^t\}$ the corresponding flow. Given two sections $\Sigma_0$ and $\Sigma_1$ of $\{\Phi_H^t\}$, we can define
the associated \textit{Poincar\'e map} $R_H\co \Sigma_0\to\Sigma_1$ by taking intersection of $\Sigma_1$ and orbits from $\Sigma_0$.  We take
$\Sigma_0$ and $\Sigma_1$ to be sections such that
the Poincar\'e map $R_H\co\Sigma_0\to\Sigma_1$
is a diffeomorphism.
Let $y_0\in\Sigma_0$ and $x_0$ be a point on the trajectory $\{\Phi_H^t(y_0)\}$ between $\Sigma_0$ and $\Sigma_1$. 

For any $h\in\mathbb{R}$ and $i\in\{0,1\}$, denote by $\Sigma_i^h:=\Sigma_i\cap H^{-1}(h)$.  $\Sigma_i^h$ is a smooth $(2n-2)$-dimensional symplectic manifold whose symplectic form is given by restriction of $\omega$. Moreover, the Poincar\'e map $R_H\co\Sigma_0\to\Sigma_1$, $R_H$ sends $\Sigma_0^h$ to $\Sigma_1^h$. We denote by $R_H^h$ the restriction $R_H|_{\Sigma_0^h}$.  It is a symplectomorphism with respect to the symplectic form defined above.

Let $\widetilde R$ be a perturbation of $R_H$ with the same
properties as of $R_H$, namely,
\begin{enumerate}
 \item $\widetilde R\co\Sigma_0\to\Sigma_1$ is a diffeomorphism;
 \item $\widetilde R$ preserves $H$, that is, $H\circ\widetilde R = H$ on $\Sigma_0$. Equivalently, $\widetilde R(\Sigma_0^h)=\Sigma_1^h$ for every $h\in\R$;
 \item the restriction of $\widetilde R$ to each $\Sigma_0^h$
is a symplectomorphism.
\end{enumerate}

Our goal in this section is to prove the following proposition.

\begin{proposition} \label{perturb all levels}


For any neighborhood $U$ of $x_0$, there exists a neighborhood $V$ of $y_0$ in $\Sigma_0$ such that for any $C^{\infty}$-small perturbation $\widetilde R$ of $R$ satisfying $(1)$-$(3)$ above and $\widetilde R=R_H$ outside $V$, there exists   a $C^{\infty}$-small perturbation $\widetilde H$ of $H$ with $\widetilde H=H$ outside $U$, and $\widetilde R = R_{\widetilde H}$.
\end{proposition}

We will prove Proposition \ref{perturb all levels} in Section \ref{s: pf of all level}. Before that we firstly realize $\widetilde R$ as a Poincar\'e map only on one level set $H^{-1}(h)$:

\begin{proposition} \label{perturb one level}
Let  $h=H(x_0)$. For any neighborhood $U$ of $x_0$ in $\Omega$
there exists a neighborhood $V^h$ of $y_0$ in $\Sigma^h_0$
such that for any $C^{\infty}$-small symplectic perturbation $\widetilde R^h$ of $R^h_H$  and $\widetilde R^h=R^h_H$ outside $V^h$, there exists a $C^{\infty}$-small perturbation $\widetilde H$ of $H$ with $\widetilde H=H$ on $H^{-1}(h)\setminus U$
and $\widetilde R^h = R^h_{\widetilde H}$.
\end{proposition}

In order to prove Proposition \ref{perturb one level} (and Lemma \ref{l:perturb section} later), we need the following generalization of Darboux's theorem 
\cite[Chapter~I, Theorem~17.2]{LM}.

\begin{theorem}[Carath\'{e}odory-Jacobi-Lie]\label{CJL}
Let $(\Omega^{2n}, \omega_0)$ be a symplectic manifold.  \textcolor{black}{Any family $\{p_1,...,p_k\}$ ($k\leq n$) of pairwise Poisson commutative and algebraically independent differentiable functions defined near $x\in \Omega$ can be completed to a symplectic coordinates $\{p_1,...,p_k, p_{k+1},...,p_{n},q_{1},...,q_{n}\}$ near $x$.}

\end{theorem}

\subsection{Outline of the proof of Proposition \ref{perturb one level}}

The proof of Propositions 
\ref{perturb one level} 
is divided into a number of steps.

\medskip
\textbf{Step 1. }\textcolor{black}{\textit{Reduction: }}
By Theorem \ref{CJL},  it suffices to prove the proposition
in this canonical setup where $\Omega=\R^{2n} = \{ (\mathbf{q},\mathbf{p}) : \mathbf{q},\mathbf{p}\in\R^n \} $, $\omega=d\mathbf{q}\wedge d\mathbf{p}$, $H= p_n$, $x_0=(\mathbf{0},\mathbf{0})$,  $\Sigma_0= \{ q_n=-1 \} $ and
$\Sigma_1= \{  q_n=1 \} $. 

\medskip
\textbf{Step 2.  }\textcolor{black}{\textit{Uniqueness of symplectomorphism determined by the targeting foliation: }
Let $\R^n_h:=\{\mathbf{p}\in\R^n: p_n=h\}.$ $\Sigma_0^h$ is foliated by Lagrangian submanifolds $\{A_{\mathbf{\widehat p}}\}_{\mathbf{\widehat p}\in\R_h^n}$ defined by $ A_{\mathbf{\widehat p}} := \{  \mathbf{p}=\mathbf{\widehat p} \} .$ Any symplectomorphism maps $\{A_{\mathbf{\widehat p}}\}_{\mathbf{\widehat p}\in\R_h^n}$ into a Lagrangian foliation of $\Sigma_1^h$.  On the other hand, the resulting Lagrangian foliation of $\Sigma_1^h$ determines the symplectomorphism by the following lemma.}

%
%

\begin{lemma}\label{lh:split}
Let $R^h_1,R^h_2\co\Sigma^h_0\to\Sigma^h_1$ be symplectomorphisms such that $R^h_1=R^h_2$ outside a compact subset of~$\Sigma^h_0$.
Suppose that
$R^h_1(A_{\mathbf{\widehat p}})=R^h_2(A_{\mathbf{\widehat p}})$
for every $\mathbf{\widehat p}\in \R_h^n$.
Then $R^h_1=R^h_2$.
\end{lemma}


{\bf Step 3.  }\textcolor{black}{\textit{Construction of Hamiltonian via Lagrangian foliation on the phase space:}
The phase space $\Omega=\R^{2n}$ is foliated by Lagrangian affine subspaces
$L_{\mathbf{\widehat p}} := \{ (\mathbf q, \mathbf{\widehat p}) : q\in\R^n \}$ as $\mathbf{\widehat p}$ ranges over $\R^{n}$.  Assume $\{\widetilde L_{\mathbf{\widehat p}}\}_{p\in\R_h^n}$ is a Lagrangian foliation such that every leaf $\widetilde L_{\mathbf{\widehat p}}$ connects $A_{\mathbf{\widehat p}}$ with its image $\widetilde R^h(A_{\mathbf{\widehat p}})$,  namely,}
\be\label{e:L-left}
 \widetilde L_{\mathbf{\widehat p}} \cap\Sigma^h_0 = A_{\mathbf{\widehat p}} \text{ and
 }\widetilde L_{\mathbf{\widehat p}} \cap\Sigma^h_1 = \widetilde R^h(A_{\mathbf{\widehat p}})\ee
for all $\mathbf{\widehat p}\in\R_h^n$.  \textcolor{black}{Module some apparent geometric restriction on $\{\widetilde L_{\mathbf{\widehat p}}\}_{p\in\R_h^n}$,  we can construct a Hamiltonian so that its Poincar\'e map is $\widetilde R^h$:}


\begin{lemma}\label{l:leaves}
Let $\{\widetilde L_{\mathbf{\widehat p}}\}_{\mathbf{\widehat p}\in\R_h^n}$ be a foliation
 by Lagrangian submanifolds satisfying \eqref{e:L-left}.
\textcolor{black}{Let $\widetilde H$ be a smooth function such that 
$\widetilde H|_{\widetilde L_{\mathbf{\widehat p}}} = h $ for all $\mathbf{\widehat p}\in\R_h^n$. Suppose $\{\widetilde L_{\mathbf{\widehat p}}\}_{p\in\R_h^n}$ satisfies some ``geometric condition", } then $R^h_{\widetilde H} = \widetilde R^h$
and
$ \widetilde H=H $
on $H^{-1}(h)\setminus U$.

\end{lemma}


{\bf Step 4.}
\textcolor{black}{\textit{Construction of $\{\widetilde L_{\mathbf{\widehat p}}\}_{p\in\R_h^n}$ satisfying Lemma \ref{l:leaves}}}. We firstly identify $\R^{2n}$ with the cotangent bundle $T^*\R^n$ using $q_i$'s are spatial coordinates and $p_i$'s are coordinates in the fibers of the cotangent bundle. 
Then we construct the desired leaves $\widetilde L_{\mathbf{\widehat p}}$ as graphs of some closed 1-forms
$\widetilde\alpha=\widetilde\alpha_{\mathbf{\widehat p}}$ on $\R^n$.  \textcolor{black}{Here we use the fact that $\widetilde R^h$ is $C^{\infty}$ close to $R^h_H$. We also manage to glue closed 1-forms derived from  $A_{\mathbf{\widehat p}}$ and $\widetilde R^h(A_{\mathbf{\widehat p}})$ so that the resulting 1-form is still closed.}

\medskip

{\bf Step 5.} \textcolor{black}{\textit{Closeness of $\widetilde{H}$ and $H$.}} Show that $\widetilde{H}$ can be chosen to be $C^{\infty}$ close to $H$ \textcolor{black}{by carefully checking the construction of the 1-forms in Step 4}.  This way we finish the proof of Proposition \ref{perturb one level}.

\subsection{Proof and application of Proposition \ref{perturb all levels}}\label{s: pf of all level}
Proposition \ref{perturb all levels} can be proved  by applying Proposition \ref{perturb one level}  to all $h\in\R$ and to the corresponding restrictions $\widetilde R|_{\Sigma_0^h}$ in place of $\widetilde R^h$. Let $\{\widetilde H^h\}_{h\in\mathbb{R}}$ be the family of the Hamiltonians generated in Proposition \ref{perturb one level}, then the desired $\widetilde H$ in Proposition \ref{perturb all levels} is constructed as $\widetilde H^{-1}(h)=(\widetilde H^h)^{-1}(h)$ for every $h\in\mathbb{R}$. 

We also apply Proposition \ref{perturb all levels} to prove the following fact which is known in folklore but for which the authors could not find a reference.  \textcolor{black}{It is used in the next section.}

\begin{proposition}\label{p:family}
Let $\varphi_0\co D^{2n}\to D^{2n}$, $n\ge 1$, be a symplectomorphism $C^\infty$-close to the identity
and coinciding with the identity near the boundary. Then there exists a smooth family
of symplectomorphisms $\{\varphi_t\}_{t\in[0,1]}$ of $D^{2n}$
fixing a neighborhood of the boundary
and such that $\varphi_t=\varphi_0$ for all $t\in[0, \frac13]$, $\varphi_t=\operatorname{id}$ for all $t\in[\frac23,1]$,
and the family $\{\varphi_t\}$ is $C^\infty$-close to the trivial family (of identity maps).
\end{proposition}

\section{Ideas of the  Proof of Theorem \ref{thm2}}\label{sec: proof}
\textcolor{black}{Let $(\mathbf q,\mathbf p)=(q_1,...,q_n, p_1,$ $..,p_n)$ be action-angle coordinates near $\mathcal T\subseteq \Omega$ with $\mathbf p=\mathbf 0$ on $\mathcal T$. By perturbing $H$ if necessary, we may assume $H(\mathbf 0)=0$, $X_H$ is colinear to $\partial/\partial p_n$ on $\mathcal T$ and the system is KAM-nondegenerate at $\mathcal T$.}

Pick a point $y_0\in\mathcal T$ and choose a small section
$\Sigma_0$ through $y_0$.  \textcolor{black}{Denote $\Sigma_1:=R(\Sigma_0)$, where $R$ is the Poincar\'e return map to $\Sigma_0$. }
The sections $\Sigma_0,  \Sigma_1$ are naturally identified with open sets in $\T^{n-1}\times D$ and parametrized by coordinates $(\bbq,\mathbf p)$
where $\bbq=(q_1,\dots,q_{n-1})$ and $\mathbf p=(p_1,\dots,p_n)$.
The explicit formula of $R$ is given by
\be\label{e:R formula}
R(\bbq,\mathbf p) = \left(q_1+\tfrac{\partial H/\partial p_1}{\partial H/\partial p_n}(\mathbf p),\dots,q_{n-1}+\tfrac{\partial H/\partial p_{n-1}}{\partial H/\partial p_n}(\mathbf p), \mathbf p \right) .
\ee
Note that the origin of $\Sigma_0$ is a fixed point of $R$.  By Abramov's formula \cite{Ab2}, in order to get Theorem \ref{thm2}, it suffices to construct a perturbed Hamiltonian so that the Poincar\'{e} map has positive volume entropy and is entropy non-expansive.  With Proposition \ref{perturb all levels}, it boils down to the following lemma:

\begin{lemma}\label{l:perturb section}
There exists a diffeomorphism $\widetilde R\co\Sigma_0 \to \Sigma_1$ arbitrarily
close to $R$ in $C^\infty$ and such that $\widetilde R=R$
outside an arbitrarily small neighborhood of the origin and the following
conditions are satisfied:
\begin{enumerate}
 \item For every $h\in\R$ near 0, 
 $\widetilde R$ maps the level set $\Sigma_0^h$ to $\Sigma_1^h$
 and is  symplectic.
 \item There is a small $\widetilde R$-invariant neighborhood of the origin and
 the restriction of $\widetilde R$ to this neighborhood has positive volume entropy.
 Moreover, $\widetilde R$ is entropy non-expansive.
\end{enumerate}
\end{lemma}

\begin{proof}[Idea of Proof]
The key ideas in the proof of Lemma \ref{l:perturb section} is divided into the following 4 steps:

\textbf{Step 1.} We first show that the restriction of $R$ on $\Sigma_0^h$ is conjugate to the following map near the fixed points: 
\be\label{e:G_h formula} 
  G_h(\bbq,\bbp) = (\bbq + \nabla f_h(\bbp), \bbp)
\ee
where $f_h$ is a function on a neighborhood of $\bbz\in\R^{n-1}$.  $G_h$ can be regarded as the time-1 map of the Hamiltonian flow 
$\Phi^t_{F_h}$ with the Hamiltonian $F_h$ given by
\be\label{e:F_h}
 F_h(\bbq,\bbp):=f_h(\bbp)
\ee
It suffices to perturb $F_h$ and define $\widetilde G_h$
as the time-1 map of the flow defined by the perturbed Hamiltonian.

\textbf{Step 2.} The KAM-nondegenerate condition at $\mathcal T$
implies that $\bbz\in\R^{n-1}$ is a nondegenerate critical point of~$f_0$.
Moreover, for all $h$ near $0$, the function $f_h$ has a nondegenerate
critical point $\bbc(h) = (c_1(h),\dots,c_{n-1}(h))$ depending smoothly on $h$
with $\bbc(0)=\bbz$. By the Morse-Bott Lemma, there exists a coordinate chart 
$\bar{\textbf{P}}^h=(P_1,\dots,P_{n-1})$ in $O_p$, depending smoothly on $h$,
such that $\bar{\textbf{P}}^h$ vanishes at $\bbc(h)$ and 
\be\label{e:morse}
f_h=f_h(\bbc(h))+P_1^2+\cdots+P_k^2-P_{k+1}^2-\cdots-P_{n-1}^2
\ee
in a neighborhood of $\bbc(h)$. Then we use Theorem \ref{CJL} to extend this collection
of functions to a symplectic coordinate system 
$(\bar{\mathbf Q}^h,\bar{\mathbf P}^h)=(Q_1,\dots,Q_{n-1},P_1,\dots,P_{n-1})$
in a neighborhood of the point $(\bbq,\bbp)=(\bbz,\bbc(h))$ for any $h$ close to 0.

\textbf{Step 3.} For any $\ep, \delta>0$, define a perturbed Hamiltonian $F_{h,\ep,\delta}$ on the $\delta$-neighborhood of $(\bbz,\bbc(h))$ by
\be\label{e:F_h,ep,delta}
F_{h,\ep,\delta}
:=F_h + \ep\,\xi(\bar{|\mathbf P}^h|)\,\xi(|\bar{\mathbf Q}^h|)\,(Q_1^2+\cdots+Q_k^2-Q_{k+1}^2-\cdots-Q_{n-1}^2),
\ee
where $\xi$ is a smooth function on $[0,1]$ with $\xi\equiv 1$ on $[0,\delta/2]$ and $\xi\equiv 0$ on $[\delta,1]$.  The perturbed Hamiltonian is periodic near $(\bbz,\bbc(h))$ with  period $2\pi/\sqrt{\ep}$.  Fix a small $h_0$, we choose $\ep(h)$ so that:
\begin{itemize}
\item When $|h|<h_0/3$, $\ep(h)$ is a constant such that  $N=2\pi/\sqrt{\ep(h)}$ is an integer.
\item When $|h|>2h_0/3$,  $\ep(h)=0$.
\end{itemize} 
Define $G_{h,\ep, \delta}:=\Phi^1_{F_{h,\ep(h), \delta}} $, and  choose a closed disc $B$ such that the sets
$B, G_{h,\ep, \delta}(B),$ $\dots, G_{h,\ep, \delta}^{N-1}(B)$ 
are disjoint for all $|h|<h_0/3$.

\textbf{Step 4.} By a result in \cite{BT2}, there exists a symplectomorphism $\theta\co \R^{2n-2}\to \R^{2n-2}$
arbitrarily $C^\infty$-close to the identity such that $\theta(B)=B$ and $\theta|_B$ has positive volume entropy.  By Proposition \ref{p:family}, there exists a smooth family $\{\theta_h\}$ that is $C^\infty$-close to $\operatorname{id}$ and such that $\theta_h=\theta$ for $|h|<h_0/3$,
and $\theta_h=\operatorname{id}$ for $|h|>2h_0/3$.

The union of compositions $\widetilde G_h = G_{h,\ep, \delta} \circ \theta_h$ on each $\Sigma_0^h$ is the desired perturbation $\tilde{R}$ of $R$. Since all $\widetilde G_h$ has positive volume entropy on $\Sigma_0^h$ when $|h|<h_0/3$, by the Abramov-Rokhlin entropy Formula \cite{AbR}, $\tilde{R}$ also has positive volume entropy.  This finishes the proof of Lemma \ref{l:perturb section}.
\end{proof}

{\it Remark}.  In Step 3,  when $n=2$, a result by Gelfreich-Turaev \cite{GT10}  shows that one can create a periodic spot near any elliptic periodic point by a $C^{\infty}$-small perturbation.  It is still unknown whether this can be generalized to higher dimensions.  Here we use special properties of integrable systems to create periodic spots in any dimensions.





\end{document}